\documentclass[12pt,a4paper,reqno]{amsart}
\usepackage{amsmath,amsthm,amssymb}

\usepackage{color}
\usepackage[dvips]{graphicx}
\theoremstyle{plain}
\newtheorem{theorem}{Theorem}
\newtheorem{lemma}{Lemma}

\newtheorem{corollary}{Corollary}
\newtheorem{definition}{Definition}
\theoremstyle{definition}
\newtheorem{example}{Example}

\usepackage [right=2.5cm, left=2.5cm, top=2.5cm, bottom=2.5cm]{geometry}

\begin{document}
\title[On open weakly $2$-semiconvex sets in three-dimensional Euclidean space]{On open weakly $2$-semiconvex sets with nonempty $2$-nonsemiconvexity-point sets in three-dimensional Euclidean space}
\author{Tetiana M. Osipchuk}
\address{Institute of Mathematics of NASU, 01004 Ukraine, Kiev-4\\
3, Tereschenkivska st.}
\email{osipchuk.tania@gmail.com}
\begin{abstract}
We study general topological properties of open weakly $2$-semiconvex sets in real three-dimensional Euclidean space $\mathbb{R}^3$. An open set in $\mathbb{R}^3$ is called \textit{weakly $2$-semiconvex} if through every boundary point of the set there passes a closed half-plane that does not intersect the set. The largest subset of the complement of a given set in $\mathbb{R}^3$ such that every closed half-plane intersecting this subset also intersects the given set is called the \textit{$2$-nonsemiconvexity-point set} corresponding to the given one. We prove that the nonempty $2$-nonsemiconvexity-point set corresponding to an open weakly $2$-semiconvex set is open
and weakly $2$-semiconvex. We also construct a weakly $2$-semiconvex domain whose $2$-nonsemiconvexity-point set is an unbounded domain that can be either convex or nonconvex.
\end{abstract}

\keywords{weakly $2$-semiconvex sets, $2$-nonsemiconvexity-point sets, three-dimensional Euclidean space.}
\subjclass{52A30}

\begingroup
\def\uppercasenonmath#1{}
\maketitle
\endgroup
\section{Introduction}
Weakly $m$-semiconvex sets, where $m=1,2,\ldots,n-1$, in the $n$-dimensional real Euclidean space $\mathbb{R}^n$, $n\geqslant2$, can be viewed as a generalization of linearly accessible domains from the theory of univalent functions. A domain in the plane is linearly accessible if its complement is a union of rays \cite{She}. The notion was introduced  by Yurii Zelinskii \cite{Zel2}, \cite{Zel3}.

An open subset $E\subset\mathbb{R}^n$, $n\geqslant2$, is called \emph{weakly $m$-semiconvex}, $m=1,2,\ldots,n-1$, if, for every point $x\in \partial E$, there exists
a closed $m$-dimensional half-plane passing through $x$ and not intersecting $E$. For a subset $E\subset\mathbb{R}^n$, $n\geqslant2$, the largest subset $E_m^{\diamondsuit}\subset \mathbb{R}^n\setminus E$, where $m=1,2,\ldots,n-1$, such that every closed $m$-dimensional half-plane intersecting $E_m^{\diamondsuit}$ intersects $E$ as well is called the \textit{$m$-nonsemiconvexity-point set} corresponding to $E$. Notice that a $1$-dimensional half-plane in $\mathbb{R}^n$, $n\geqslant2$, is a ray, and a $2$-dimensional half-plane in $\mathbb{R}^3$ is just a half-plane.

The theory of weakly $m$-semiconvex sets is relatively young, and the results obtained so far primarily concern weakly $1$-semiconvex sets. In particular, it is known that, for every open weakly $1$-semiconvex subset $E\subset\mathbb{R}^n$, $n\geqslant2$, with $E_1^{\diamondsuit}\ne\varnothing$, the set $E_1^{\diamondsuit}$ is open and weakly $1$-semiconvex, \cite{Osi2024}, \cite{Osi2025}. Moreover, an open weakly $1$-semiconvex subset $G\subset \mathbb{R}^2$ with nonempty $1$-nonsemiconvexity-point set $G_1^{\diamondsuit}$ has at least three connected components \cite{Dak}. Furthermore, each connected component of $G_1^{\diamondsuit}$ is bounded and convex, and every connected subset of its boundary consisting entirely of smooth points is either a line segment or a single point  \cite{Osi2024}. By contrast, there exists an open weakly $1$-semiconvex subset $E\subset\mathbb{R}^n$, $n>2$, such that the set $E_1^{\diamondsuit}$ is connected and nonconvex and may be either bounded or unbounded  \cite{Osi2025}. Examples of weakly $m$-semiconvex domains in $\mathbb{R}^n$, $n\geqslant 3$, with nonempty $m$-nonsemiconvexity-point sets are also known for $1\leqslant m<n-1$ \cite{Osi2022}.

Thus, the problem of finding a lower bound for the number of connected components of weakly $(n-1)$-semiconvex sets in $\mathbb{R}^n$, $n \geqslant 3$, has remained open, as has the question regarding the topological properties of $m$-nonsemiconvexity-point sets corresponding to weakly $m$-semiconvex sets in $\mathbb{R}^n$, where $n\geqslant 3$ and $1< m\leqslant n-1$.

In this paper, we investigate analogous properties of weakly $2$-semiconvex sets in $\mathbb{R}^3$. In particular, we prove that, for every open weakly $2$-semiconvex subset $E\subset\mathbb{R}^3$ with $E_2^{\diamondsuit}\ne\varnothing$, the set $E_2^{\diamondsuit}$ is open and weakly $2$-semiconvex. We also construct a weakly $2$-semiconvex domain $E\subset\mathbb{R}^3$ for which $E_2^{\diamondsuit}$ is an unbounded domain that can be either convex or nonconvex.

\section{Definitions and Preliminaries}

An $m$-dimensional affine subspace of $\mathbb{R}^n$, $n\geqslant1$, where  $m=0,1,2,\ldots,n-1$, is called an {\bf\emph{$m$-dimensional plane}}.

An {\bf\emph{ $m$-dimensional half-plane}} in $\mathbb{R}^n$, $n\geqslant2$, where  $m=1,2,\ldots,n-1$, is one of the two closed parts into which an $m$-dimensional plane is divided by any of its $(m-1)$-dimensional planes.

We use the following standard notations. For a subset $A\subset\mathbb{R}^n$, let $\overline{A}$
be its closure, $\mathrm{Int}\, A$ be its interior, and $\partial
A=\overline{A}\setminus\mathrm{Int}\, A$ be its boundary. 

\begin{definition}[cf. \cite{Zel2}]\label{def1} An open subset $E\subset\mathbb{R}^n$, $n\geqslant2$, is called {\bf\emph{weakly $m$-semiconvex}}, $m=1,2,\ldots,n-1$, if, for every point $x\in \partial E$, there exists
an $m$-dimensional half-plane $H$ such that $x\in H$ and $H\cap E=\varnothing$. The class of weakly $m$-semiconvex sets in $\mathbb{R}^n$ is denoted by $\mathbf{WS^n_m}$.
\end{definition}



\begin{definition}[cf. \cite{Osipchuk2019}]\label{def2}
A point $x\in \mathbb{R}^n\setminus E$ is called an {\bf\emph{$m$-nonsemiconvexity}} point of a subset $E\subset\mathbb{R}^n$ if every $m$-dimensional half-plane $H$ containing $x$ intersects $E$. The set of all $m$-nonsemiconvexity points of $E$ is called the {\bf\emph{$m$-nonsemiconvexity-point set}} corresponding to $E$ and is denoted by $E_m^{\diamondsuit}$.
\end{definition}

In Definitions \ref{def1} and \ref{def2},  the condition $x\in H$ can be replaced by $x\in \partial H$, without changing the respective definitions.

The class of weakly $m$-semiconvex sets with nonempty $m$-nonsemiconvexity-point sets is denoted by $\mathbf{WS^n_m}\setminus \mathbf{S^n_m}$. The sets belonging to $\mathbf{WS^n_m}\setminus \mathbf{S^n_m}$ have the following properties.
\begin{lemma}[Dakhil \cite{Dak}, Theorem 2.4.1]
An open set in $\mathbb{R}^2$ belonging to the class $\mathbf{WS^2_{1}}\setminus \mathbf{S^2_{1}}$ has at least three connected components.
\end{lemma}

Interestingly, the number of connected  components of a set belonging to the class $\mathbf{WS^2_{1}}\setminus \mathbf{S^2_{1}}$ also depends on the smoothness of its boundary.
We say that a component $G$ of an open bounded set in $\mathbb{R}^2$ has {\bf\emph{smooth boundary}} if $\partial G$ is the image of a $C^1$-embedding of the unit circle.
We say that an open bounded set in $\mathbb{R}^2$ has smooth boundary if each of its components has smooth boundary.

\begin{lemma}[Osipchuk \cite{Osipchuk2019}, Theorem 1.8]
Suppose that an open bounded subset $E\subset \mathbb{R}^2$ with smooth boundary belongs to the class $\mathbf{WS^2_1}\setminus \mathbf{S^2_1}$. Then $E$
has at least four connected components.
\end{lemma}

Examples of open sets belonging to the class $\mathbf{WS^2_{1}}\setminus \mathbf{S^2_{1}}$ with three connected components, as well as examples of such sets with smooth boundaries and at least four connected components, can be found in  \cite[Examples~2.1--2.4]{Osi2022}.


\begin{lemma}[Osipchuk \cite{Osi2022}, Theorem 6]
For every $m=1,2,\ldots, n-2$, there exists a domain in $\mathbb{R}^n$, $n\geqslant 3$, belonging to the class $\mathbf{WS^n_m}\setminus \mathbf{S^n_m}$.
\end{lemma}

The following results describe properties of the $1$-nonsemiconvexity-point sets corresponding to weakly $1$-semiconvex sets in $\mathbb{R}^n$, $n\geqslant2$.

\begin{lemma}[Osipchuk \cite{Osi2024}, Theorem 1.2 (a); \cite{Osi2025}, Theorem 2.2]
Let $E\subset \mathbb{R}^n$, $n\geqslant 2$, be an open subset belonging to the class $\mathbf{WS^n_1}\setminus \mathbf{S^n_1}$. Then $E^{\diamondsuit}_1$ is open.
\end{lemma}
\begin{lemma}[Osipchuk \cite{Osi2024}, Theorem 1.2 (e); \cite{Osi2025}, Corollary 2.3]
Let $E\subset \mathbb{R}^n$, $n\geqslant 2$, be an open subset belonging to the class $\mathbf{WS^n_1}\setminus \mathbf{S^n_1}$. Then $E^{\diamondsuit}_1\in\mathbf{WS^n_1}$.
\end{lemma}
\begin{lemma}[Osipchuk \cite{Osi2024}, Theorem 1.2 (b)-(d)]
Let $E\subset \mathbb{R}^2$ be an open subset belonging to the class $\mathbf{WS^2_1}\setminus \mathbf{S^2_1}$. Then
\begin{enumerate}
\item[$\bullet$] the components of $E^{\diamondsuit}_1$ are convex and bounded;
\item[$\bullet$] every connected subset of $\partial E^{\diamondsuit}_1$ consisting entirely of smooth points is either a line segment or a single point.
\end{enumerate}
\end{lemma}
\begin{lemma}[Osipchuk \cite{Osi2025}, Theorem 2.11]
There exists an open subset $E\subset \mathbb{R}^n$, $n\geqslant 3$, belonging to $\mathbf{WS^n_1}\setminus \mathbf{S^n_1}$, such that the set $E_1^{\diamondsuit}$ is bounded (or unbounded), connected, and nonconvex.
\end{lemma}

\section{Main results}
Given two points $x,y\in\mathbb{R}^3$, we denote by $xy$ the open line segment joining $x$ and $y$, by $\|x-y\|$ the length of this segment, and by $\langle x,y\rangle$ the Euclidean inner product of $x$ and $y$. Let also
$$S^{2}:=\{z\in \mathbb{R}^3: \|z\|=1\};$$
$$U(y,r):=\{x\in\mathbb{R}^3:\|x-y\|<r\},\quad y\in\mathbb{R}^3,\quad r>0;$$
$$
\gamma_a(y):=\{\alpha a+y: \alpha\in(-\infty,+\infty)\},
$$
$$\eta_a(y):=\{\alpha a +y: \alpha\in[0,+\infty)\},$$
$$L(y,a):=\{x\in \mathbb{R}^3: \langle a,x-y\rangle=0\},\quad
y\in\mathbb{R}^3,\quad a \in S^{2};$$
\begin{multline*}
H(y,a,b):=\{x\in \mathbb{R}^3: x-y=\alpha a+\beta b,\,\,\alpha\in(-\infty,+\infty),\,\,\beta\in[0,+\infty)\},\\ y\in\mathbb{R}^3,\,\, a,b \in S^{2},\,\,\langle a,b\rangle=0.
\end{multline*}
Let also $O$ be the origin of $\mathbb{R}^3$, and set
$$
H(a,b):=H(O,a,b).
$$
Observe that $H(a,b)=H(-a,b)$ as subsets of $\mathbb R^3$. Nevertheless, we consider the ordered pairs $(a,b)$ and $(-a,b)$ as distinct elements of the parameter space
$$
\mathcal{H}:=\{(a,b)\in S^{2}\times S^{2}: \langle a,b\rangle=0\}
$$
equipped with the subspace topology inherited from $\mathbb{R}^6$.

Define
$$
SO(3):=\{\Gamma\in\mathbb R^{3\times3}: \Gamma^T\Gamma=\mathrm{I}_3,\ \det \Gamma=1\},
$$
where $\Gamma^T$ denotes the transpose of $\Gamma$, $\mathrm{I}_3$ is the $3\times3$ identity matrix, and $\det \Gamma$ is the determinant of $\Gamma$.

\begin{lemma}\label{lem1}
For an open ball $U(p,r)\subset \mathbb{R}^3$, the subset
$$V :=\left\{(a,b) \in \mathcal{H}: H(a,b)\cap U(p,r) \ne\varnothing\right\}$$
is open in $\mathcal{H}$.
\end{lemma}
\begin{proof}[Proof.] Define $D:\mathcal H\to[0,+\infty)$ as the distance function from $p$ to $H(a,b)$, i.e.,
$$
D(a,b):=\inf\limits_{x\in H(a,b)}\|p-x\|.
$$
Then the condition
$$
H(a,b)\cap U(p,r) \ne\varnothing
$$
is equivalent to
$$
D(a,b)<r.
$$
Consequently,
$$
V=\{(a,b)\in\mathcal H: D(a,b)<r\}.
$$
Show that $D$ is continuous on $\mathcal H$. Since $a$, $b$ and $a\times b$ form an orthonormal basis of $\mathbb R^3$, we can write
$$
p=\langle p,a\rangle a+\langle p,b\rangle b+
\langle p,a\times b\rangle(a\times b).
$$
Every point $x\in H(a,b)$ has the form $x=\alpha a+\beta b$, where $\alpha\in\mathbb R$ and $\beta\geqslant0$. Hence
$$
\|p-x\|^2
=(\langle p,a\rangle-\alpha)^2+(\langle p,b\rangle-\beta)^2
+\langle p,a\times b\rangle^2.
$$
The minimum is attained at
$$
\alpha=\langle p,a\rangle,\qquad
\beta=\max\{\langle p,b\rangle,0\}.
$$
Therefore,
$$
D(a,b)^2
=\langle p,a\times b\rangle^2
+\min\{\langle p,b\rangle,0\}^2.
$$
Geometrically, this means that the distance from $p$ to $H(a,b)$ equals the distance from $p$ to the plane $L(O,a\times b)\supset H(a,b)$ if the projection of $p$ on $L(O,a\times b)$ lies in $H(a,b)$, and equals the distance from $p$ to the boundary of $H(a,b)$ if the projection of $p$ on $L(O,a\times b)$ is not in $H(a,b)$.
Since
$$
\min\{t,0\}=\frac{t-|t|}{2},
\qquad t\in\mathbb R,
$$
and the inner product, the cross product, and the absolute value are continuous, the function $D^2$ is continuous on $\mathcal H$. Consequently, $D=\sqrt{D^2}$ is also continuous on $\mathcal H$.
Finally,
$$
V=D^{-1}((-\infty,r)).
$$
Since $(-\infty,r)$ is open in $\mathbb R$ and $D$ is continuous, $V$ is open in $\mathcal H$.
\end{proof}
\begin{lemma}\label{lem2}
Let $E\subset \mathbb{R}^3$ be an open subset such that $E_2^{\diamondsuit}\ne\varnothing$. Let $y\in E_2^{\diamondsuit}$. Then for every $(a,b) \in \mathcal H$, there exists an open ball $U(x(y,a,b),d(y))\subset E$, where $x(y,a,b)\in H(y,a,b)\cap E$ and $d(y)>0$ depends only on $y$, not on $(a,b)$.
\end{lemma}
\begin{proof}[Proof.]
Fix an arbitrary point $y\in E_2^{\diamondsuit}$.
Since $y\in E_2^\diamondsuit$, every half-plane $H(y,a,b)$, $(a,b)\in\mathcal H$, intersects $E$. For each $(a,b)\in\mathcal H$, choose a point $p(a,b)\in H(y,a,b)\cap E$.
Since $E$ is open, there exist open balls $U_{a,b}:=U\left(p(a,b),r(a,b) \right)$, $(a,b) \in \mathcal H$, such that $\overline{U_{a,b}}\subset E$.

Let
$$V_{a,b}:=\left\{ (c,g) \in \mathcal H: H(y,c,g)\cap U_{a,b} \ne\varnothing\right\},\quad (a,b) \in \mathcal H.$$

By Lemma~\ref{lem1}, the sets $V_{a,b}$, $(a,b) \in \mathcal H$, are open, moreover, $(a,b) \in V_{a,b}$. Then
$$
\bigcup\limits_{(a,b) \in \mathcal H} V_{a,b}
$$
is a cover of $\mathcal H$.

Consider the mapping $\Phi:\mathcal{H}\to SO(3)$, defined by
$$
\Phi(a,b)=(a,b,a\times b),
$$
where the vectors on the right-hand side are the columns of the corresponding matrix, is a homeomorphism. Indeed, $\Phi$ is bijective, since every matrix in $SO(3)$ is uniquely determined by its first two columns, while its third column is their cross product. Both $\Phi$ and its inverse are continuous, as their coordinate functions are continuous.

Since $\mathcal{H}$ is homeomorphic to $SO(3)$, and $SO(3)$ is compact; see \cite{Lee}, $\mathcal{H}$ is compact as well.
Then, by the Heine-Borel theorem, there exists a subcover
$$
\bigcup\limits_{j\in M}
V_{a_j,b_j},\quad (a_j,b_j) \in \mathcal H,\quad j\in M,\quad M\,\,\mbox{is finite},
$$
of $\mathcal H$. Consequently, any half-plane $H(y,a,b)$, $(a,b) \in \mathcal H$, intersects $\bigcup\limits_{j\in M}\overline{U_{a_j,b_j}}$.

Let $E_i$, $i\in N\subseteq \mathbb{N}$, be the components of $E$. Then for any $j\in M$ there exists $i(j)\in N$ such that $\overline{U_{a_j,b_j}}\subset E_{i(j)}$.

Consider the distance functions
$$d_i(x):=\inf\limits_{x^0\in \partial E_i}\|x-x^0\|,\quad x\in E_i,\quad i\in N.$$
They are continuous in the domains $E_i$, $i\in N$. Then their restrictions to the compacts $\overline{U_{a_j,b_j}}$ attain their minimum values $d_{j}>0$ on these compacts, i.e.,
$$
d_{j}:=\min\limits_{x\in \overline{U_{a_j,b_j}}}d_{i(j)}(x),\quad j\in M.
$$
Since $M$ is finite, there exists
$$
d:=\min\limits_{j\in M}d_j>0.
$$

\begin{figure}[h]
\centering
\includegraphics[width=15 cm]{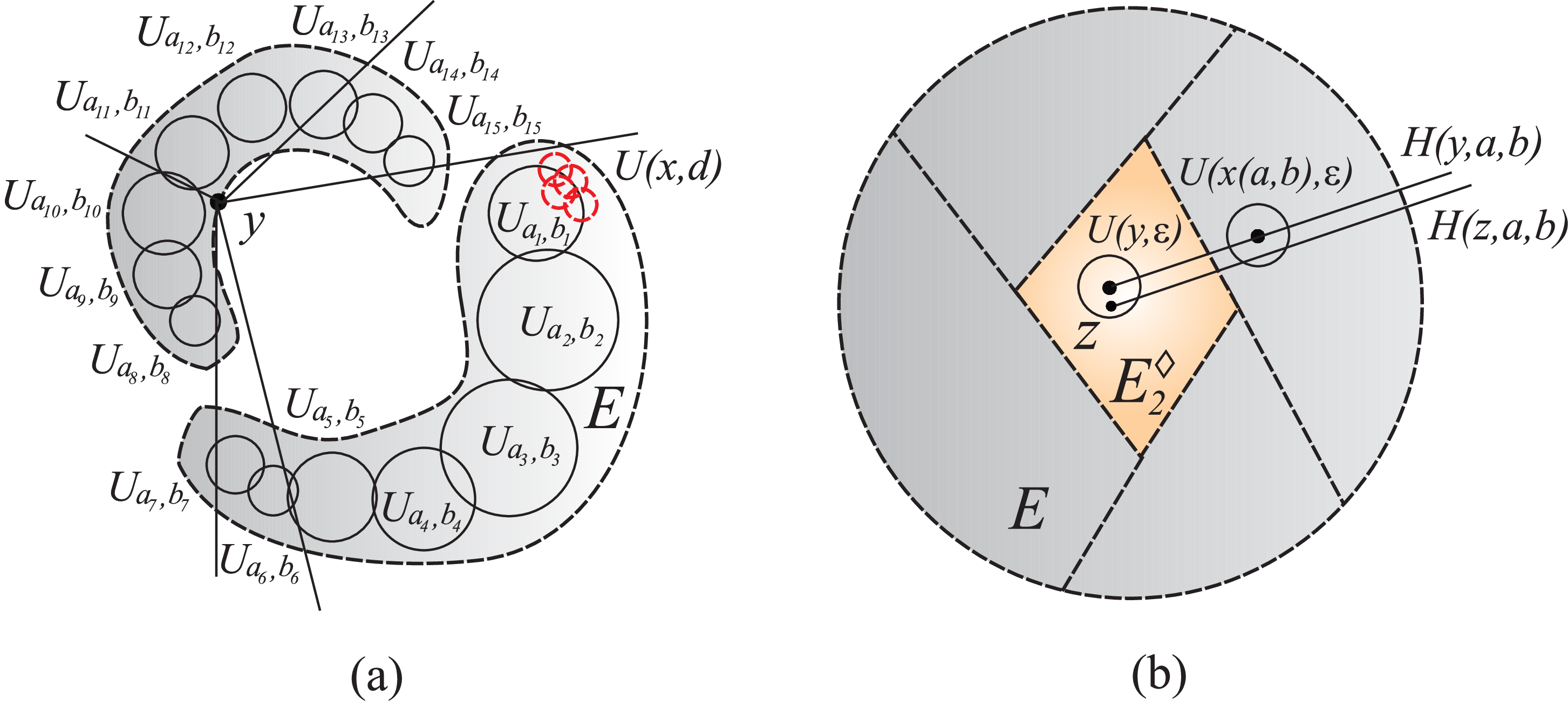}
\caption{}\label{Fig4_1}
\end{figure}

By the construction of the finite subcover, for every $(a,b)\in\mathcal H$, there exists $j\in M$ such that $H(y,a,b)\cap U_{a_j,b_j}\ne\varnothing$.
Choose a point
$x(y,a,b)\in H(y,a,b)\cap U_{a_j,b_j}$.
Since $x(y,a,b)\in\overline{U_{a_j,b_j}}$, we have
$$
U(x(y,a,b),d)\subset E.
$$
Thus, the required assertion follows.
\end{proof}
\begin{theorem}\label{theor1}
Let an open subset $E\subset \mathbb{R}^3$ belong to the class $\mathbf{WS^3_2}\setminus \mathbf{S^3_2}$. Then $E_2^{\diamondsuit}$ is open.
\end{theorem}
\begin{proof}[Proof.] Fix an arbitrary point $y\in E_2^{\diamondsuit}$ and show that it is an inner point of $E_2^{\diamondsuit}$.

Since $E$ is weakly $2$-semiconvex, it follows that $y\not\in \partial E$. Then there exists a number $\varepsilon_1>0$ such that $U(y,\varepsilon_1)\subset \mathbb{R}^3\setminus \overline{E}$.

By Lemma~\ref{lem2}, for the fixed $y$ and any $(a,b) \in \mathcal H$ there exist a point $x(a,b) \in H(y,a,b)\cap E$ and a constant $d>0$ such that $U(x(a,b),d)\subset E$.

Let $\varepsilon:=\min\{\varepsilon_1,d\}$. Consider the neighborhood $U(y, \varepsilon)$ of the point $y$. Let $z\in U(y, \varepsilon)$, and $H(z,a,b)$, $(a,b) \in \mathcal H$, be an arbitrary half-plane passing through $z$. Draw the half-plane $H(y,a,b)$ parallel to the half-plane $H(z,a,b)$. Since $U(x(a,b), \varepsilon)\subseteq U(x(a,b), d)\subset E$ for the point $x(a,b)$ correspondent to $H(y,a,b)$, it follows that $H(z,a,b)\cap U(x(a,b), \varepsilon)\ne\varnothing$ and, therefore, $H(z,a,b)\cap E\ne\varnothing$ for any $(a,b) \in \mathcal H$; see Figure~\ref{Fig4_1}\,(b). Thus, $z$ is a $2$-nonsemiconvexity point of $E$. Since $z$ is arbitrary, it implies that all points of $U(y, \varepsilon)$ are $2$-nonsemiconvexity points of $E$. Hence, $y$ is an inner point of $E_2^{\diamondsuit}$.
\end{proof}
\begin{corollary}\label{corol1}
Let an open subset $E\subset \mathbb{R}^3$ belong to the class $\mathbf{WS^3_2}\setminus \mathbf{S^3_2}$. Then $E_2^{\diamondsuit}$ is weakly $2$-semiconvex.
\end{corollary}
\begin{proof}[Proof] Since the set $E_2^{\diamondsuit}$  is open, it follows that $\partial E_2^{\diamondsuit}\cap E_2^{\diamondsuit}=\varnothing$. Moreover, $\partial E_2^{\diamondsuit}\cap E=\varnothing$, since $E$ is open and $E_2^{\diamondsuit}\subset \mathbb{R}^3\setminus E$. Then, through every point $y\in\partial E_2^{\diamondsuit}$, there passes a half-plane $H(y,a,b)$, $(a,b) \in \mathcal H$, that does not intersect $E$. In this case, $H(y,a,b)\cap E_2^{\diamondsuit}=\varnothing$ by Definition~\ref{def2}. Thus, $E_2^{\diamondsuit}$ is weakly $2$-semiconvex.
\end{proof}

\begin{example}
Let us construct the example of a domain $E\in \mathbf{WS^3_2}\setminus \mathbf{S^3_2}$ such that the set $E^{\diamondsuit}_2$ is an unbounded domain that can be either convex or nonconvex.

We identify $\mathbb{R}^2$ with the horizontal plane $\mathbb{R}^2\times \{0\}$. Let $T\subset \mathbb{R}^2$ be an open curvilinear triangle with vertices $y^j$, $j=1,2,3$, such that the sides $y^1y^2$, $y^1y^3$ are line segments. Denote by $\gamma(y^1y^2)$ and $\gamma(y^1y^3)$ the straight lines containing the segments $y^1y^2$ and $y^1y^3$, respectively. Let $A$ be the open angle bounded by $\gamma(y^1y^2)$ and $\gamma(y^1y^3)$ and containing $T$. Let also $y^2y^3$ be a curve without common points with straight lines $\gamma(y^1y^2)$, $\gamma(y^1y^3)$ except the endpoints $y^2$ and $y^3$, and such that every ray $\eta_a(y^1)$ emanating from $y^1$ and contained in $A$ intersects the curve $y^2y^3$ exactly at one point; see Figure~\ref{Fig1}.

Consider the sets
$$
P:=T\times(-\infty,0),\quad G:=(\overline{A\setminus T})\times (-\infty,+\infty),
$$
and the half-planes
$$
H^2:=\gamma(y^1y^2)\times(-\infty,0],\quad H^3:=\gamma(y^1y^3)\times(-\infty,0].
$$
Let also $H^1$ be the closed half-plane in $\mathbb{R}^2$ containing $T$ and whose boundary passes through $y^1$ parallel to the straight line $\gamma(y^2y^3)$ passing through $y^2$ and $y^3$.
Let
$$
H:=\bigcup\limits_{k=1}^3H^k.
$$
Consider the set
$$
E:=\mathbb{R}^3\setminus \left(P\cup G \cup H\right).
$$
It is a domain by the construction.
\begin{figure}[h]
\centering
\includegraphics[width=8 cm]{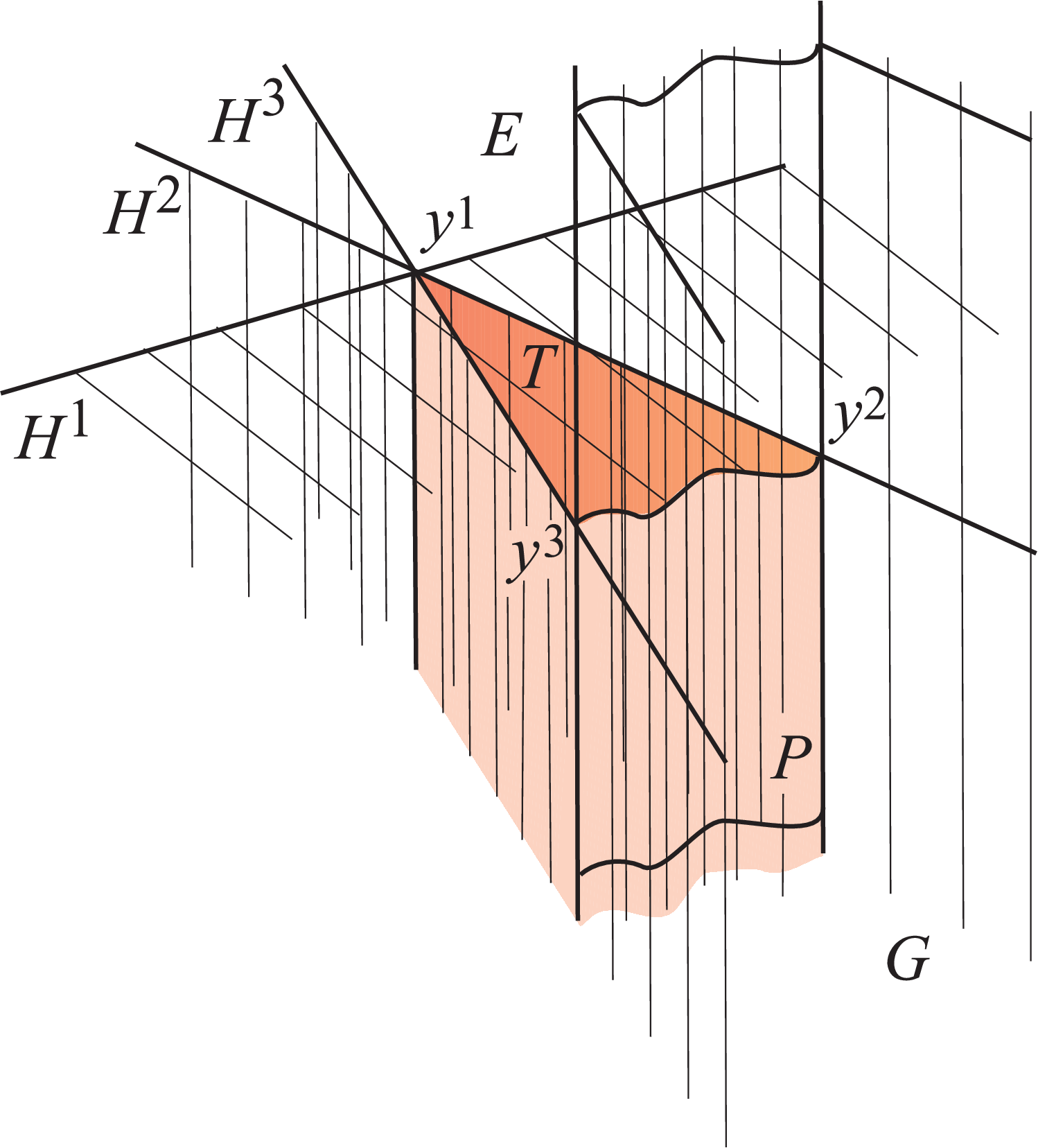}
\caption{}\label{Fig1}
\end{figure}
Show that $E\in \mathbf{WS^3_2}$. Considering the general properties of set boundaries and the fact that $\partial P\subset \partial G\cup H$ by the construction,
we obtain that
$$
\partial E\subset \partial P\cup\partial G \cup H=\partial G\cup H.
$$
It is obvious that for any point $x\in H$ there exists a half-plane in $H$ passing through $x$ and, therefore, not intersecting $E$.

If $x\in G$, then consider the orthogonal projection $x'$ of $x$ onto $\mathbb{R}^2$. The half-plane $\eta_a(x')\times (-\infty,+\infty)$, where the direction of $a$ is defined by the ray emanating from $y^1$ and passing through $x'$, is completely contained in $G$.  This follows from the condition of the unique intersection of any ray $\eta_a(y^1)$ and the curve $y^2y^3$. Therefore, $\eta_a(x')\times (-\infty,+\infty)\cap E=\varnothing$. Since $G$ is closed, we obtain that $\partial G\cap E^{\diamondsuit}_2=\varnothing$. Thus, $E\in \mathbf{WS^3_2}$.

Since
$$
E^{\diamondsuit}_2\subset\mathbb{R}^3\setminus E=P\cup G \cup H,
$$
it also follows that $E^{\diamondsuit}_2\subset P$.
Show that $P\subset E^{\diamondsuit}_2$. First, prove that $\gamma_a(x)\cap E\ne \varnothing$ for any $a\in S^2$, $x\in P$. Consider two cases:
\begin{enumerate}
\item[1.] If $a=(0,0,-1)$ or $a=(0,0,1)$, the straight line $\gamma_a(x)$ is vertical. In this case,
$$
\gamma_a(x)\cap \left[T\times (0,+\infty)\right]\ne\varnothing, \quad x\in P.
$$
Since $\left[T\times (0,+\infty)\right]\subset E$, it implies that
$$
\gamma_a(x)\cap E\ne\varnothing \quad x\in P.
$$
\item[2.] If $a\in S^2\setminus\left(\{(0,0,-1)\}\cup\{(0,0,1)\}\right)$, the straight line $\gamma_a(x)$ is not vertical. Consider the set
$$
\widetilde{E}:= \mathbb{R}^3\setminus \left(\left(A\cup\gamma(y^1y^2)\cup\gamma(y^1y^3)\right)\times \mathbb R\cup H^1\right).
$$
Then, for any point $x\in P\subset  A\times\mathbb R$,
$$
\gamma_a(x)\cap \widetilde{E}\ne\varnothing.
$$
Indeed, the orthogonal projection of $\gamma_a(x)$ onto $\mathbb{R}^2$  is the straight line passing through the orthogonal projection of $x$. Since $x\in A\times\mathbb R$, this line passes through a point of $A$ and, being a complete straight line, intersects
$$\mathbb{R}^2\setminus \left[A\cup\gamma(y^1y^2)\cup\gamma(y^1y^3)\right].$$
Moreover, $\gamma_a(x)\not\subset H^1$, since $x\not \in H^1$.

Since $\widetilde{E}\subset E$, it follows that
$$
\gamma_a(x)\cap E\ne\varnothing,\quad x\in P.
$$
\end{enumerate}

Since every half-plane $H(x,a,b)$ contains the straight line $\gamma_a(x)$, we conclude that
$$
H(x,a,b)\cap E\ne\varnothing
$$
for every $x\in P$ and every $(a,b)\in\mathcal H$. Thus, $P\subset E_2^\diamondsuit$.

Therefore, we proved that $E^{\diamondsuit}_2= P$. By choosing $T$ to be convex, nonconvex, or a straight-sided triangle, we obtain a convex, nonconvex, or even prism-shaped set $P$, respectively.
\end{example}

\section*{Declarations}
This work was supported by the National Research Foundation of Ukraine in the framework of the research project 2025.07/0014 on the topic "Modern problems of Mathematical Analysis and Geometric Function Theory". The author declares no potential conflict of interest with respect to the research, authorship and publication of this article. All necessary data are included into the paper.

\end{document}